\documentclass[english,12pt]{article}
\usepackage[totalwidth=430pt,totalheight=595pt]{geometry}
\usepackage{babel,amsmath,amssymb,amsbsy,amsfonts,latexsym,amsthm,epsf,array,colordvi,xy,fancyhdr,mathrsfs}
\input xy
\xyoption{all}

\newtheorem{teor}{Theorem}[section]
\newtheorem{lemma}[teor]{Lemma}
\newtheorem{prop}[teor]{Proposition}

\newcommand{\ac}{\`{a} }

\newcommand{\graor}{\mathbb{G}_{3}}
\newcommand{\gra}{\mathbb{G}_{3}}
\newcommand{\grasc}[1]{\mathbb{G}_{2}(\C^{\,#1})}
\newcommand{\grasr}[1]{\mathbb{G}_{4}(\R^{\,#1})}
\newcommand{\la}{\mathfrak{g}}

\newcommand{\sud}{\mathfrak{su}(2)}

\newcommand{\soq}{\mathfrak{so}(4)}
\newcommand{\spu}{\mathfrak{sp}(1)}
\newcommand{\spd}{\mathfrak{sp}(2)}

\newcommand{\C}{\mathbb{C}}
\newcommand{\Ha}{\mathbb{H}}
\newcommand{\R}{\mathbb{R}}
\newcommand{\Z}{\mathbb{Z}}

\newcommand{\Pro}{\mathbb{P}}

\newcommand{\taut}{\mathbf{V}}

\newcommand{\eu}{\R^{n}}

\newcommand{\grak}{\mathbb{G}_{k}}

\newcommand{\twistor}{\mathcal{Z}}
\newcommand{\swannbundle}{\mathcal{U}}

\newcommand{\si}[1]{\Sigma^{#1}}
\newcommand{\conn}[1]{\nabla^{#1}}
\newcommand{\mlift}{\hat{\Phi}}

\newcommand{\hh}{\hat{h}}
\newcommand{\hbeta}{\hat{\beta}}
\newcommand{\sym}[1]{\{#1\}}

\newcommand{\graskn}[2]{\mathbb{G}_{#1}(\R^{\,#2})}

\newcommand{\grad}{\mathrm{grad}}

\newcommand{\vspan}{\mathrm{span}}

\newcommand{\vhom}{\mathrm{Hom}}
\newcommand{\vend}{\mathrm{End}}

\newcommand{\twistqk}{\mathscr{D}}
\newcommand{\killing}{\mathscr{K}}

\newcommand{\hol}{\mathrm{Hol}}
\newcommand{\Tr}{\mathrm{Tr}}
\newcommand{\gop}{\boldsymbol{\gamma}}
\newcommand{\conf}{\boldsymbol{\lambda}}
\newcommand{\hpsi}{\hat{\Psi}}

\newcommand{\nqd}{\mathcal{N}_{\Ha}}
\newcommand{\autnatural}[1]{(#1)^{\natural}}
\newcommand{\uf}{\theta}
\newcommand{\sect}{s}
\newcommand{\ufd}{\zeta}
\newcommand{\summ}{C}

\title{\bfseries Latent Quaternionic Geometry}
\author{\textsc{Andrea Gambioli}}
\date{10th April 2006}

\begin{document} 
\maketitle

\abstract{\footnotesize In this article we discuss the interaction between the geometry of a quaternion-K\"{a}hler
manifold $M$ and that of the Grassmannian $\graor(\la)$ of oriented $3$-dimensional subspaces of a compact Lie algebra
$\la$. This interplay is described mainly through the moment mapping induced by the action of a group $G$ of quaternionic
isometries on $M$.  We give an alternative expression for the endomorphisms $I_{1},I_{2},I_{3}$, both in terms of the
holonomy representation of $M$ and the structure of the Grassmannian's tangent space. A correspondence between the solutions
of respective twistor-type equations on $M$ and $\graor(\la)$ is provided.  \vspace{4mm}

}

\section{Introduction}

Let $G$ be a compact Lie group acting by quaternionic isometries on a qua\-ter\-nion-K\"ahler (QK) manifold $M$. In this case a
Killing vector field $X$ satisfies the condition $L_X\Omega=0$, where $\Omega$ is the parallel $4$-form of the QK
structure. Recall that the fibre of the standard rank $3$ vector bundle over $M$ (whose complexification is often written
$S^{2}H$) is isomorphic to $\spu$, and is spanned by a basis of endomorphisms $I_{1}, I_{2}, I_{3}$ satisfying the
quaternionic relations \begin{equation*} I_{i}^{2}=-\hbox{Id}\quad \text{and}\quad I_{i}I_{j}=\epsilon^{ijk} I_{k}
\end{equation*} with $\epsilon ^{ijk}$ the sign of the permutation.

We denote by $\mu$ the moment map for the $G$ action, and by $\mu_{A}$ the section of $S^{2}H$ obtained by the contraction
of $\mu$ with $A\in\la$ through the metric induced by the Killing form. It satisfies the equation
\begin{equation}\label{twisteqnqk}
d \mu_{A}\,=i(\tilde{A})\,\Omega\,,
\end{equation}
where $\tilde{A}$ is the Killing vector field generated by $A$ (see \cite{gal}, \cite{galaw1}). Another way of describing the sections
coming from the moment map is expressed by the formula
\begin{equation*}
\mu_{A}=\pi_{S^{2}H}(\nabla\tilde{A})
\end{equation*}
up to a constant. The moment map $\mu$ is $G$-equivariant with respect to the given action of $G$ on $M$ and of the adjoint
representation of $G$ on $\la$: it can be used to construct the $G$-equivariant morphism \[\Psi:M_{0}\xymatrix{\ar[r]&}
\graor(\la),\] where $M_{0}$ is an appropriate subset of $M$. The morphism $\Psi$ was introduced by Swann (\cite{swann91},
\cite{swann98}), who studied the unstable manifolds for the gradient flow of an appropriate functonal $\psi$ on this type of
Grassmannians, proving that they admit a QK structure; we will use the map $\Psi$ in order to relate in various ways the
geometry of QK manifolds to that of Grassmannians of type $\graor(\la)$.

 In Section \ref{operatorsongrassmannians}, we introduce the natural first-order differential operator $D$ on the
tautological rank $k$ vector bundle over a Grassmannian $\graskn{k}{n}$, which annihilates projections of constant
sections. Indeed, we show that all solutions of $D$ arise in this way (Theorem~2.2). This illustrates a well-known
technique, whereby solutions of an overdetermined differential operator may be interpreted as parallel sections of some
connection on a larger bundle (\cite{bceg}).  Although quaternionic geometry and Lie algebras are not yet involved, we aim to show
that $D$ is completely analogous to the more complicated \emph{twistor operator} $\twistqk$ on a QK manifold.

 In Section \ref{thetwotwistorequations}, we recall the definition of $\twistqk$ on sections of the vector bundle $S^2H$,
and explain that it is satisfied by the moment sections $\mu_A$ defined above. We then prove that under suitable hypotheses
the map $\Psi$ can be used to relate elements in $\ker\twistqk$ with those in $\ker D$ where $D$ now acts on the
tautological rank 3 vector bundle $V$ over $\graor(\la)$.

 Whilst the tangent space to $\graor(\la)$ at $V$ is given by \begin{equation}\label{VVperp} T_{V}\graor(\la)\cong
V\otimes V^{\perp},\end{equation} the complexified tangent space to $M$ has the form $H\otimes_\C E$, reflecting the
representation of the holonomy group $Sp(1)Sp(n)$. Part of our problem is to reconcile the roles of the ``auxiliary'' vector
bundles $V$ and $H$ with respective fibres $\R^3$ and $\C^2$. In Section \ref{thespnsp1structure} we give an alternative
description of the imaginary quaternion endomorphisms $I_{i}$ over a point $x\in M$ in terms of $Sp(1)$ representations of
a subgroup $Sp(1)$ diagonally embedded in $Sp(1)Sp(n)$.

In Section \ref{thecoincidencetheorem} we state our main results: we show that it is possible to push forward the
endomorphisms $I_{1},I_{2},I_{3}$ so that they can be described as endomorphisms of the subspace $\Psi_{*}T_{x}M$ of
(\ref{VVperp}), where $V=\Psi(x)$. In other words, if $Z=\sum_{i=1}^{3} v_{i}\otimes p_{i}$ belongs to $\Psi_{*}T_{x}M$,
then we can write \begin{equation*} I_{k}Z=\sum_{i=1}^{3} v_{i}\otimes q_{i} \end{equation*} and we shall explicitly
determinine the $q_{i}$s in terms of $v_{i},\,p_{i}$ and $I_{k}$. This is accomplished in Proposition 5.4, itself a
geometric counterpart to the re\-pre\-sen\-ta\-tion-theoretic Proposition~4.1.

 Finally, in Section \ref{examplesandapplications} we apply the theory to the case of an $Sp(1)\times Sp(1)$ action on
$\Ha\Pro^{1}$ and to other compatible examples. We describe some natural real 4-dimensional subspaces of (\ref{VVperp})
which correpond to quaternionic lines in $T_xM$, and are tangent to quaternion projective lines in Wolf spaces.


\section{Operators on Grassmannians}\label{operatorsongrassmannians}

Consider an $n$-dimensional real vector space $\eu$ equipped with an inner product $\langle\,,\,\rangle$; we can construct
the Grassmannian of oriented $k$-planes $\grak(\eu)$, whose tangent space at a $k$-plane $V$ can be identified with the
linear space \begin{equation*} \vhom(V,V^{\perp})\cong V^{*}\otimes V^{\perp}\,; \end{equation*} in fact if
$v_{1},\cdots,v_{k}$ is an orthonormal (ON) basis for $V$ and $w_{1},\cdots,w_{n-k}$ for $V^{\perp}$, then each homomorphism
$T_{ij}$ defined as $T_{ij}(v_{k})=\delta_{k}^{i}w_{j}$, corresponds to an independent tangent direction; more explicitly,
the curve \begin{equation}\label{tangij} \alpha_{ij}(r):=\vspan\{v_{1},\cdots,(\cos{r})v_{i}+(\sin{r})w_{j},\cdots,v_{k}\}
\end{equation} satisfies $\alpha_{ij}(0)=V$ and $\alpha_{ij}'(0)=T_{ij}$. The presence of a metric on $V$,\,induced from the
ambient space $\eu$, will allow us to write $V\otimes V^{\perp}$,\,using contraction via the metric for the isomorphism
$V\cong V^{*}$.

We will be interested in studying differential operators and sections of vector bundles on $\graskn{k}{n}$, so we start
by describing some induced objects. Given the metric, we have the splitting of the trivial bundle $\grak(\eu)\times \eu$ in
two subbundles:\,the tautological one $\mathbf{V}$ and its orthogonal complement: \begin{equation*} \xymatrix{
\mathbf{V}\oplus\mathbf{V}^{\perp}\ar ^-{\cong} [r]\ar[ddr]_{p} &\grak(\eu)\times\eu\ar[dd]^{p'}&\\ & \\ & \grak(\eu) }\,.
\end{equation*} The presence of this metric also allows us to define connections on these two subbundles merely by
composing $d$ with the two projections $\pi$ and $\pi^{\perp}$. For instance \begin{equation*}
\nabla^{\mathbf{V}}s=\pi \,d \,s\,, \end{equation*} where $s\in\Gamma(\taut)$ and $d$ is the derivation in $\eu$. To prove
that this is a connection let $a$ be a function, and note that \begin{align*} \nabla^{\taut}(as)=\pi
d(as)&=\pi\left((da)s+a(ds)\right) \notag \\ &=(da)s+a \pi (ds) \notag\\ &=(da)s+a\nabla^{\taut}s \end{align*} as
required.\,Moreover this connection is compatible with the metric induced on the fibres of $\taut$ by their ambient space
$\eu$: in fact if $s,t\in\Gamma(\taut)$ and $X\in T_{V}\grak(\eu)$ we have \begin{align*} X\langle s\,,\,t\rangle =\langle
Xs\,,\,t\rangle +\langle s\,,\,Xt\rangle =&\langle\pi Xs\,,\,t\rangle +\langle s\,,\,\pi Xt\rangle \notag \\
=&\langle\nabla^{\mathbf{V}}_{X}s\,,\,t\rangle +\langle s\,,\,\nabla^{\mathbf{V}}_{X}t\rangle \,.  \end{align*}

On the other hand we obtain the corresponding second fundamental form by projecting in the opposite 
way:
\begin{equation*}
\Gamma(\taut)\xymatrix{\ar[r]&}\Gamma(T^{*}\grak(\eu)\otimes\taut^{\perp})
\end{equation*}
which sends $s$ to $\pi^{\perp}ds$; analogously $\mathit{II}^{\perp}$ sends $s\in\Gamma(\taut^{\perp})$ to $\pi ds$. 
Both $\mathit{II}$ and $\mathit{II}^{\perp}$ are tensors. In fact, if for example $s\in\Gamma(\taut^{\perp})$
and $a$ is a function, we get
\begin{equation*}
\pi d (as)=\pi (d(a)s+ad(s))=\pi ad(s)=a \pi ds
\end{equation*}
so that we can think to $\mathit{II}^{\perp}$ as a section of the bundle 
\begin{equation*}
\vhom\left(\taut^{\perp}\,,\,T^{*}\grak(\eu)\otimes\taut\right)\cong \taut^{\perp}\otimes\left(T^{*}\grak(\eu)\otimes\taut\right)
\end{equation*}
(identifying $\taut^{\perp}\cong(\taut^{\perp})^{*}$ as usual). It turns out that this section determines an
immersion of $\taut^{\perp}$ as a subbundle of $T^{*}\grak(\eu)\otimes\taut$; we shall return to this question later in 
the Section.
\vspace{2mm}

We use the standard connections and tensors previously introduced in order to construct new differential
operators on the tautological bundle $\taut$ and on its orthogonal complement $\taut^{\perp}$. 
First of all, given an element $A\in\eu$ we can associate to it two sections of the bundles $\taut$ and $\taut^{\perp}$ 
just using the projections: $s_{A}=\pi A$ and $s^{\perp}_{A}=\pi^{\perp}A$ with $A=s_{A}+s^{\perp}_{A}$; as $A$ is constant,
\begin{equation*}
0=dA=ds_{A}+ds^{\perp}_{A}
\end{equation*}
so that
\begin{equation*}
ds_{A}=-ds^{\perp}_{A}\,;
\end{equation*}
in the language already deployed
\begin{equation*}
\nabla^{\taut}s_{A}=\pi ds_{A}=-\pi d s^{\perp}_{A}=-\mathit{II}^{\perp}s^{\perp}_{A}\,.
\end{equation*}
These equations imply that 
\begin{equation}\label{costsecimm}
d \,s_{A}=-\mathit{II}^{\perp}s^{\perp}_{A}+\mathit{II}s_{A}\,.
\end{equation}

For convenience we will combine the homomorphisms $\mathit{II}$ and $\mathit{II}^{\perp}$ to act upon any
$\eu$-valued function on $\graor(\eu)$, giving a mapping
\begin{equation*}
i:\xymatrix{C^{\infty}(\graor(\eu),\eu)\ar[r]&\Gamma(T^{*}\otimes\eu})
\end{equation*}
defined by
\begin{equation}\label{immerstaut}
i(S)=\mathit{II}(\pi S)-\mathit{II}^{\perp}(\pi^{\perp}S)\,.
\end{equation}
in a way which is consistent with equation (\ref{costsecimm}). Thus we have
\begin{equation}\label{costsec1}
ds_{A}=i(A)
\end{equation}
and
\begin{equation}\label{costsec2}
ds_{A}^{\perp}=-i(A)\,.
\end{equation}\medbreak

The image of $\mathit{II}^{\perp}$ corresponds to elements of the type
\begin{equation}\label{vperpingr}
\sum_{i=1}^{k} \lambda\,y\otimes v_{i}\otimes v_{i}
\end{equation}
with $y \in \taut^{\perp}$ and $\lambda\in\R$; this can be shown with the following argument: let us consider the 
decomposition as $SO(k)\times SO(n-k)$ modules of the involved bundles
\begin{equation}\label{decomptaut}
\taut^{\perp}\otimes\taut\otimes\taut\cong\taut^{\perp}\otimes \R+\taut^{\perp}\otimes(\taut\otimes\taut)_{0}
\end{equation}
where $(\taut\otimes\taut)_{0}$ is the tracefree part of the tensor product; Schur's Lemma guarantees that the second summand 
cannot contain any submodule isomorphic to $\taut^{\perp}$, so the first summand consists of the unique submodule of
this type in the right side term of (\ref{decomptaut}). Therefore, as expression (\ref{vperpingr}) provides an 
$SO(k)\times SO(n-k)$-equivariant copy of $\taut^{\perp}$ inside this bundle, it must coincide with 
$\mathit{II}^{\perp}(\taut^{\perp})$. The same argument shows that 
\begin{equation*}
\mathit{II}(u)=\sum_{i=1}^{n-k}\lambda\, u \otimes w_{i}\otimes w_{i} 
\end{equation*} 
with $u\in \taut,\,\lambda\in \R$. We want now to be more precise about these statements, and calculate explicitly the value
of $\lambda$. This is done in the next proposition (in which tensor product symbols are omitted).

\begin{prop}\label{costantprojsec}
Let $A\in\eu$ so that $A=u+y$ with $u\in V$ and $y\in V^{\perp}$ at the point $V$; let $v_{j}$ and $w_{i}$ denote the 
elements of ON bases of $V$ and $V^{\perp}$ at $V$;\,then 
\begin{equation}\label{sff1}
\mathit{II}(u)=\sum_{j}{u\,w_{j}w_{j}}
\end{equation}
and
\begin{equation}\label{sff2}
\mathit{II}^{\perp}(y)=-\sum_{i} {y\,v_{i}v_{i}}\,.
\end{equation}
\end{prop} 
\noindent\textit{Proof.} We differentiate the section $s_{A}$ along the curve $\alpha_{ij}(t)$ passing through $V$ and with 
tangent vector $v_{i}w_{j}$ as in (\ref{tangij}); let $u=\sum_{i=1}^{k} a_{i}v_{i}$ and $y=\sum_{j=1}^{n-k}b_{j}w_{j}$; then
\begin{align*}
s_{A}(\alpha_{ij})(t)&=a_{1}v_{1}+\cdots+\langle A,\cos{r}\,v_{i}+\sin{r}\,w_{j}  \rangle(\cos{r}\,v_{i}+\sin{r}\,w_{j})+\cdots+v_{k}\notag\\
&=a_{1}v_{1}+\cdots+(a_{i}\cos{r}+b_{j}\sin{r})(\cos{r}\,v_{i}+\sin{r}\,w_{j})+\cdots+v_{k}\notag
\end{align*}  
so that 
\begin{equation*}
\frac{d}{dr}s_{A}(\alpha_{ij})(r)_{|_{r=0}}=d\,s_{A}\cdot v_{i}w_{j}=b_{j}v_{i}+a_{i}w_{j}\,;
\end{equation*}
therefore, as an $\R^{n}$-valued $1$-form,
\begin{align*}
d\,s_{A}&= \sum_{ij}b_{j}\,v_{i}v_{i}w_{j}+ a_{i}\, w_{j}v_{i}w_{j}\notag\\
&=\sum_{i} y \,v_{i}v_{i}+ \sum_{j}u \,w_{j}w_{j}\,,
\end{align*}
where the second summand belongs to $\taut \otimes \taut^{\perp}\otimes\taut^{\perp}$ and coincides with $\mathit{II(u)}$ 
as claimed. An analogous calculation for $s_{A}^{\perp}$ gives
\begin{equation*}
d\,s_{A}^{\perp}=-\sum_{i}{yv_{i}v_{i}}-\sum_{j}{uw_{j}w_{j}} 
\end{equation*}
as expected from equation (\ref{costsec2}).\,$\blacksquare$
\vspace{2mm}

 \noindent\textit{Observation.} The opposite signs in (\ref{sff1}) and (\ref{sff2}) are consistent with the equation
$$
0=d\langle s_{A}\,,\, s_{A}^{\perp}\rangle|_{V}=\langle \mathit{II}(u)\,,\,y\rangle+\langle u\,,\,\mathit{II}^{\perp}(y)\rangle 
$$
which expresses the fact that $\mathit{II}$ and $\mathit{II}^{\perp}$ are adjoint linear operators.
\vspace{2mm}

Proposition \ref{costantprojsec} shows that $\nabla^{\taut}s_{A}$ is of the form seen in (\ref{vperpingr}),\,or
alternatively that if we call $\pi_{2}$ the projection on the second summand in the decomposition (\ref{decomptaut}) and define 
$D\equiv \pi_{2}\circ\nabla^{\taut}$, the section $s_{A}$ satisfies the {\em twistor-type} equation
\begin{equation}\label{twisteq}
D\,s_{A}=0\,.
\end{equation}
Symmetrically we can define another operator $D^{\perp}$ such that 
\begin{equation}\label{twisteqperp}
D^{\perp}\,s^{\perp}_{A}=0\,.
\end{equation}

Let us choose an orthonormal basis $e_{1},...,e_{n}$ of $\eu$, every section $S$ of the flat bundle $\grak(\eu)\times\eu$
is nothing else than an $n$-tuple of functions 
\begin{equation*}
f_{j}:\grak(\eu)\xymatrix{\ar[r]&}\eu
\end{equation*}
so that
\begin{equation*}
S=\sum f_{j}e_{j}\,;
\end{equation*}
applying the exterior derivative on $\eu$ (which is a connection on the flat bundle) we obtain
\begin{equation*}
dS=\sum df_{j}\otimes e_{j}
\end{equation*}
and if $1\wedge i$ denotes an element in $\vhom\left(T^{*}\otimes\eu,(\bigotimes^{2}T^{*})\otimes\eu\right)$ 
(where $T^{*}=T^{*}\grak(\eu)$ to lighten the notation) acting in the obvious way, we obtain
\begin{equation*}
1\wedge i\,(dS)=\sum df_{j}\wedge i\,(e_{j})\,;
\end{equation*}
on the other hand 
\begin{equation*}
d\sum f_{j}\,i(e_{j})=\sum df_{j}\wedge i(e_{j})+f_{j}\,di(e_{j})\,,
\end{equation*}
so if we can show that
\begin{equation*}
di(e_{j})=0\quad\forall j
\end{equation*}
we obtain the commutativity of the following diagram:
\begin{equation}\label{diagtwistinj}
\xymatrix{
     & \eu \ar[r]^-{d}\ar[d]_-{i} &T^{*}\otimes\eu\ar[d]^-{1\wedge i\quad\quad ;}&\\                  
\eu\ar[r]^-{d}&T^{*}\otimes\eu \ar[r]^-{d}&\Lambda^{2}T^{*}\otimes \eu  \\
}\;
\end{equation}
but equation (\ref{costsec1}) implies:
\begin{equation*}
di(e_{j})=dds_{e_{j}}=0\,,
\end{equation*}
because the $e_{j}$ are constant.
A consequence of Proposition \ref{costantprojsec} is that $i$ is an injective map (because $\mathit{II}$ and 
$\mathit{II}^{\perp}$ are); if we can show that also $1\wedge i$ is injective (and it happens to be in most part of cases, as 
we will see) looking at diagram (\ref{diagtwistinj}) we can deduce the following facts: if $s\in \Gamma(\taut)$
satisfies $Ds=0$, then $ds=i(s+s')$ for some $s' \in\Gamma(\taut^{\perp})$; this follows by
comparing
\begin{equation*}
ds=\nabla\,s +\mathit{II}(s)
\end{equation*}
with (\ref{immerstaut}) and noting that $\pi\,s=s$ in this case: then $s'=-(\mathit{II}^{\perp})^{-1}(\nabla\,s)$. 
Obviously $dds=0$, so $d(s+s')=0$ too, hence it is a constant element $A\in\eu$. This implies the main result of this Section:
\begin{teor}\label{teorkerntwist}
A section $s\in \Gamma(\taut)$ satisfies the twistor equation $Ds=0$ if and only if exists another section $s'\in 
\Gamma(\taut^{\perp})$ such that $s+s'=A$ is a constant section of $\eu$, provided $k>1$ and $n-k>1$.
\end{teor}

In other words sections of type $s_{A}$ are the only solutions of equation (\ref{twisteq}), under these hypotheses.
\vspace{2mm}

The missing piece to prove Theorem \ref{teorkerntwist} is injectivity of $1\wedge i$.\,To prove that we start defining another 
map:
\begin{equation*}
c:\Gamma(T^{*}\otimes\eu)\xymatrix{\ar[r]&}\Gamma(\eu)
\end{equation*}
acting as a contraction in the following way:
\begin{equation*}
c\Big(\sum_{ijk} a^{ijk} \;v_{i}w_{j}v_{k}+\sum_{lmo} b^{lmo}\; w_{l}v_{m}w_{o}\Big)
=\sum_{ij} a^{iji}\; w_{j}+\sum_{lm} b^{lml}\;v_{m}\,. 
\end{equation*}
The same map acts also on $\tau\in(\bigotimes^{q} T^{*})\otimes \eu$ in the following way: if $\tau=\tau'\otimes \uf$
with $\tau'\in\bigotimes^{q-1}T^{*}$ and $\uf\in T^{*}\otimes\eu$ then
\begin{equation*}
c(\tau)=\tau'\otimes c(\uf)
\end{equation*}
and then extending linearly.\\*
We are now in position to prove the previously stated assertion, which concludes the proof of Theorem \ref{teorkerntwist}: 
\begin{lemma}
The map $1\wedge i$ is injective,\,provided $k>1$ and $n-k>1$.
\end{lemma}
 \noindent\textit{Proof.} Given two bases $v^{i}$ of $V$ and $w^{j}$ of $V^{\perp}$ an element in $T^{*}\otimes\eu$ is described by
\begin{equation*}
\tau=\sum_{ijh} a^{ijh}\,v_{i}w_{j}v_{h}+\sum_{lmo}b^{lmo}v_{l}w_{m}w_{o} \,;
\end{equation*}
now we will prove that $c\circ 1\wedge i$ is injective,\,so that $1\wedge i$ must be.\\*
So we get
\begin{align*}
1\wedge i\;(\tau)&=\sum_{ijh\mu} a^{ijh}\,(v_{i}w_{j}\wedge v_{h}w_{\mu})w_{\mu}+\sum_{lmo\nu}b^{lmo}(v_{l}w_{m}\wedge w_{o}v_{\nu})v_{\nu}\\
               &=\sum_{ijh\mu} a^{ijh}\,(v_{i}w_{j}\otimes v_{h}w_{\mu}-v_{h}w_{\mu}\otimes v_{i}w_{j})w_{\mu}\\
               &+\sum_{lmo\nu}b^{lmo}(v_{l}w_{m}\otimes w_{o}v_{\nu}-v_{o}w_{\nu}\otimes w_{l}v_{m})v_{\nu}
\end{align*}
and applying the contraction
\begin{align*}
c(1\wedge i\;(\tau))&=\sum_{ijh\mu} a^{ijh}\,(v_{i}w_{j}\otimes v_{h}-v_{h}w_{\mu}\otimes v_{i}\delta^{j}_{\mu})\\
               &+\sum_{lmo\nu}b^{lmo}(v_{l}w_{m}\otimes w_{o}-v_{o}w_{\nu}\otimes w_{l}\delta^{m}_{\nu})\,.
\end{align*}
Now imposing that it's zero,\,we get the following couples of equations:
\begin{equation*}
\begin{cases}
(n-k)\, a^{ijh}-a^{hji}=0\\
(n-k)\, a^{hji}-a^{ijh}=0 \\
\end{cases}
\end{equation*}
and
\begin{equation*}
\begin{cases}
k\, b^{lmo}-b^{oml}=0\\
k\, b^{oml}-b^{lmo}=0 \\
\end{cases}
\end{equation*}
which imply 
\begin{equation*}
(n-k)^{2}\, a^{ijh}=a^{ijh}
\end{equation*}
and
\begin{equation*}
k^{2} \,b^{lmo}=b^{lmo}
\end{equation*}
which are absurd if $k>1$ and $n-k>1$.\;$\blacksquare$


\section{The two twistor equations}\label{thetwotwistorequations}

Let us consider a compact Lie group $G$ acting by isometries on a QK manifold $M$; then its moment map 
$\mu$ can be described locally as
\begin{equation}\label{momentmaploc}
\mu=\sum_{i=1}^{3}\omega_{i}\otimes B_{i}
\end{equation}
with $\omega_{i}$ a local orthonormal basis for $S^{2}H$ and $B_{i}$ belonging to $\la$. Suppose that
$V:=\vspan\{B_{1}, B_{2}, B_{3}\}$ is a $3$-dimensional subspace of $\la$: then $V$ is independent of the 
trivialization, as the structure group of $S^{2}H$ is $SO(3)$. Therefore we obtain a well defined map
\begin{equation*}
\Psi:M_{0}\xymatrix{\ar[r]&}\graor(\la)
\end{equation*}
where $M_{0}\subset M$ is defined as the subset where $V(x)$ is $3$-dimensional; this turns out to be an open dense
subset of the union $\bigcup S$ of $G$-orbits $S$ on $M$ such that $\dim S\geq3$ (\cite[Proposition 3.5]{swann98}). Therefore
if the dimension of the maximal $G$ orbits in $M$ is big enough, then $M_{0}$ is an open dense subset of $M$.
\vspace{3mm}

\noindent{\bfseries Assumption.} From now on we will assume that 
\begin{equation}\label{assumptionconf}
B_{i}=\conf(x) v_{i}
\end{equation}
for $v_{i}$ an orthonormal basis of $V$. This hypothesis is not excessively restrictive, in the sense that it is 
compatible with the existence of open $G_{\C}$ orbits on the twistor space $\twistor=\Pro(\swannbundle)$: in fact 
the projectivization of the complex-contact moment map $f$ induced on $\twistor$ satisfies 
\begin{equation*}
(\Pro f)(\omega_{1})=\vspan_{\C}\{B_{2}+\imath B_{3}\}\,,
\end{equation*}
and in this case this turns out to be a ray of nilpotent elements in $\la_{\C}$ (see (\cite[\S 3]{swann98}). Nilpotent elements belong to the
zero set of any invariant symmetric tensor over $\la_{\C}$, in particular with respect to the Killing form: for by Engel's
theorem their adjoint representation can be given in terms of strictly upper triangular matrices, with respect to a suitable
basis, and the product of such matrices is still strictly upper triangular and hence traceless; in other words \begin{align*}
0=\Tr\,(ad_{B_{2}+\imath B_{3}}\circ ad_{B_{2}+\imath B_{3}})&=\langle B_{2}+\imath B_{3}\,,\, B_{2}+\imath B_{3}\rangle \notag \\
&=\|B_{1}\|^{2}-\|B_{2}\|^{2}+2\imath\langle B_{2}\,,\,B_{3}\rangle\,,
\end{align*}
which implies $B_{2}\perp B_{3}$ and $\|B_{2}\|=\|B_{3}\|$, conditions that are equivalent to the assumption, permuting
cyclically the indices. Therefore condition (\ref{assumptionconf}) holds for all unstable manifolds described in 
\cite{swann98}, as in that case the twistor bundle $\twistor$ is $G_{\C}$-homogeneous. We assume throughout the 
Section that this condition holds for the moment map $\mu$. 
\vspace{2mm}

Using the map $\Psi$, we can construct on $M_0$ the pullback bundle $\Psi^{*}(\taut)$; the latter is unique up to
isomorphism of bundles (see \cite[Chap.\,I, Prop.\,2.15]{wells}). More precisely, any vector bundle $W\longrightarrow
M_{0}$ for which there exists a map of bundles $\hat{\Phi}:W\longrightarrow \mathbf{V}$ which is injective on the fibres, 
and a commutative diagram \begin{equation}\label{pullbackdiag} \xymatrix{ W\ar [r]^-{\hat{\Phi}} \ar[d]_{p^{*}_{V}} & \mathbf{V}\ar[d]^{p_{V}}\\
M_{0}\ar [r] _{\Psi} &\graor(\la),} \end{equation} is necessarily isomorphic to
$\Psi^{*}(\mathbf{V})$.
 
\begin{lemma} We have the following isomorphism of bundles on $M_{0}$: \begin{equation*} S^{2}H\cong\Psi^{*}(\mathbf{V}).
\end{equation*} \end{lemma}  \noindent\textit{Proof.} To complete the commutative diagram (\ref{pullbackdiag}), define the
morphism of bundles \begin{equation*} \xymatrix{ \hat{\Phi}:S^{2}H_{_{}} \ar[r] & \mathbf{V}} \end{equation*} by
\begin{equation*} \xymatrix{ \big(x,\omega_{i}(x)\big)\ar @{|->}[r] & \big(\vspan
\{B_{1}(x),B_{2}(x),B_{3}(x)\},B_{i}(x)\big)}\, \end{equation*} (see (\ref{momentmaploc})), extending linearly on the
fibres. This corresponds to the contraction of a vector $v\in S^{2}H_{x}$ with the $S^{2}H$ component of $\mu(x)$ using the
metric, so it does not depend on the trivialization (the structure group preserves the metric) and is injective on the
fibres by definition of $M_{0}$.\, $\blacksquare$\medbreak

We should point out that $\hat{\Phi}$ is not an isometry of Riemannian bundles in general; nevertheless under the hypotheses
discussed above, we can assume that $\hat{\Phi}$ is a conformal map of Riemannian bundles, considering $S^{2}H$ and $\taut$
to be equipped with the natural metrics coming respectively from $M$ and from $\graor(\la)$.  \vspace{2mm}

Let us now recall some well-known differential operators (the symbol $\Gamma$ denoting space of sections is omitted):\,the
{\em Dirac operator} \begin{equation*} \xymatrix{ \delta:S^{2}H \ar[r]^-{\nabla} & E\otimes H\otimes S^{2}H\, \ar @ {^{(}->}
[r] & (E\otimes {\underline H})\otimes (H\otimes {\underline H}^{*}) \ar[r] & T^{*}} \end{equation*} where the underlined
terms are contracted and $T^{*}=E\otimes H$; the {\em QK twistor operator} is defined as follows: \begin{equation*}
\xymatrix{ \twistqk:S^{2}H \ar[r]^-{\nabla} & E\otimes H\otimes S^{2}H\, \ar [r]^-{sym} & E\otimes S^{3}H}\,, \end{equation*}
where we symmetrize after covariant differentiation.\,In \cite[Lemma 6.5]{sal82}, under the assumption of nonzero scalar
curvature, Salamon proved that sections of $S^{2}H$ belonging to $\ker \,\twistqk$ are in bijection with the elements in the 
space $\killing$ of Killing vector fields preserving the QK structure; this means that if $\nu$ is in $\ker\,\twistqk$ then $\delta(\nu)$ is dual to a
Killing vector field $\tilde{A}\in \killing$,\,and on the other hand $\nu=\mu_{A}$,\,or in other words
\begin{equation}\label{twisteqquat} \twistqk \,\mu_{A}=0\, \end{equation} and all elements in $\ker\,\twistqk $ are of this
form.  \vspace{2mm}

Recall now what was discussed for Grassmannians in Section \ref{operatorsongrassmannians}: there we introduced another differential operator $D$
on the tautological bundle $\taut$ over $\graor(\la)$; the elements in its kernel were proved to be precisely the sections
$s_{A}$ obtained by projection from the trivial bundle with fibre $\la$ (see Theorem \ref{teorkerntwist}). We want to relate the kernels of
$\twistqk$ and $D$ through the map $\Psi$ induced by $\mu$; recall that the bundle homomorphism $\hat{\Phi}$ is defined up to 
a bundle automorphism of $S^{2}H$; we can for instance introduce a dilation 
\begin{equation}\label{dilationconf}
\xi(x,w)=(x,\frac{w}{\|B_{i}\|}),\,
\end{equation}
which is independent of the trivialization;\,in this way
\begin{equation*}
\hat{\Xi}(\omega_{i}):=\hat{\Phi}\circ\xi(\omega_{i})=\frac{B_{i}}{\|B_{i}\|},\,
\end{equation*}
and so an orthonormal basis is sent to another orthonormal basis: this is therefore an isometry of the two bundles compatible
with the map $\Psi$ induced by $\mu$.\,
\vspace{2mm}

We can now state the main result of this Section. Let us denote by $\killing_{\la}\subset \killing$ the subspace of Killing
vector fields induced by $\la$ and by $(\ker\,\twistqk)_{\la}$ the space of the corresponding twistor sections; then

\begin{prop}\label{proptwisteq}
There exists a lift $\hat{\Psi}$ of the map $\Psi$ such that 
\begin{equation*}
\hat{\Psi}(\mu_{A})=s_{A}\,,
\end{equation*}
inducing a bijective linear map
\begin{equation*}
(\ker\,\twistqk)_{\la}\xymatrix{\ar[r]&}\ker\,D\,.
\end{equation*}
\end{prop}
 \noindent\textit{Proof.} We are looking for a lift $\hat{\Psi}$ such that the diagram 
\begin{equation*}
\xymatrix{
S^{2}H\ar  [r]^-{\hat{\Psi}}    & \mathbf{V}\\
 M_{0}\ar  [r]  _{\Psi} \ar[u]^{\mu_{A}} &\graor(\la)\ar[u]_{s_{A}} \,.
}
\end{equation*}
commutes; recall the usual local description (\ref{momentmaploc}) of $\mu$,\,and let us define $\hat{\Psi}$ so that
\begin{equation*}
\hat{\Psi}(\omega_{i})=\frac{B_{i}}{\|B_{i}\|^{2}}\,,
\end{equation*}
obtaind by composing $\hat{\Phi}$ with the dilation $\xi^{2}$ (see (\ref{dilationconf})); this is again a lift of $\Psi$;\,consider as usual 
$\mu_{A}\in\Gamma(S^{2}H)$ satisfiying the twistor equation;\,then
\begin{align*}
\hat{\Psi}(\mu_{A})&=\hat{\Psi}\big(\sum_{i} \omega_{i}\langle B_{i}\,,\,A\rangle \big) \notag \\
&=\sum_{i}\frac{B_{i}}{\|B_{i}\|^{2}}\langle B_{i}\,,\,A\rangle\, \notag\\
&=\pi_{V}A=s_{A}\,,
\end{align*}
as required. As the lift $\hat{\Psi}$ is injective on the fibres, and as 
$$
\dim (\ker\,\twistqk)_{\la}=\dim \,\killing_{\la}=\dim\,\la=\dim\,\ker\,D
$$ 
the last assertion follows.\,$\blacksquare$
\vspace{2mm}

The situation can be summarized in diagram (\ref{diagtwistop}):
\begin{equation}\label{diagtwistop}
\xy
\xymatrix{
 &*+<13pt>[F-,:]\hbox{$A\in\la$} \POS[];[dddr]**\dir{-}& \\
& & \\
& & \\
*+<13pt>[F-,:]\hbox{$s_{A}\in \ker\,D$}\POS[];[rr]**\dir{-}\POS[];[uuur]**\dir{-}&  &*+<10pt>[F-,:]\hbox{$\mu_{A}\in (\ker\,\twistqk)_{\la}$}
}
\endxy
\end{equation}\;.
\vspace{5mm}

 \noindent\textit{Observation.} We can interpret $\mu$ as a collection of $n=\dim\,\la$ sections of $S^{2}H$:\,if $A_{i}$ are an orthonormal basis 
for $\la$ the moment map $\mu$ is completely determined by the $\mu_{A_{i}}$.\,Locally we get
\begin{equation}\label{locmmsect}
B_{i}=\sum_{j}a_{i}^{j}\,A_{j}
\end{equation}
so that 
\begin{equation*}
\mu_{A_{i}}=\sum_{j}a^{j}_{i}\,\omega_{j}\,.
\end{equation*}
\noindent For instance, if a section $\nu\in\Gamma(S^{2}H)$ is given locally by
\begin{equation*}
\nu=\sum_{i}c^{i}\omega_{i}
\end{equation*}
then 
\begin{equation*}
\hat{\Phi}(\nu)=\sum_{i}c^{i}B_{i}\,;
\end{equation*}
with respect to the basis $A_{i}$ of $\la$ the local description of the morphism $\hat{\Phi}$ is encoded in 
the $(3\times(n-3))$
matrix of the coefficients $a^{i}_{j}$ seen in $(\ref{locmmsect})$.\,\vspace{2mm}


\section{The $Sp(1)Sp(n)$ structure}\label{thespnsp1structure}

We are going now to introduce an alternative description of the endomorphisms $I_{1},I_{2},I_{3}$ in a purely 
algebraic setting, using the holonomy representation at a fixed point $x\in M$.

Let $h,\hh$ denote a unitary basis of $H$, in such a way that $\omega_{H}(h,\hh)=1$;\,with respect to this basis we have
\begin{equation}\label{invariantsymplf}
\omega_{H}=h\wedge\hh=\frac{1}{2}(h\hh-\hh h) \,.
\end{equation}
We can in terms of $h,\hh$ dtermine a basis of $S^{2}H$:
\begin{eqnarray}\label{quatbasis}
I_{1}&=&\imath(h\vee \hh) \notag\\
I_{2}&=&h^{2}+\hh^{2}\\
I_{3}&=&\imath(h^{2}-\hh^{2})\notag
\end{eqnarray}
are orthogonal of norm $\sqrt{2}$ with respect to the metric $\omega_{H}\otimes\omega_{H}$ induced on $S^{2}H$; they satisfy
the same relations of quaternions:
\begin{equation*}
I_{k}^{2}=-1\quad,\quad I_{i}I_{j}=sgn_{(ijk)}I_{k}
\end{equation*}
with $sgn_{(ijk)}$ the sign of the permutation; the composition is obtained by contracting again with $\omega_{H}$.

Consider now the case where the $Sp(1)$ representation inside $Sp(1)Sp(n)$ is such that the projection on the $Sp(n)$ factor
is nonzero: this means that the $E$ representation is nontrivial under this $Sp(1)$ action.

In this case it is significant to analyze the quaternionic action from the point of view of these new $Sp(1)$
representations. First we adopt the following notation: we have the symmetrization map $S$ acting on tensors as
\begin{equation*}
S(x_{1}\otimes\cdots\otimes x_{n})=\frac{1}{n!}\sum_{\pi^{n}}x_{\pi^{n}(1)}\otimes\cdots\otimes x_{\pi^{n}(n)}
\end{equation*}
where $\pi^{n}$ varies in the group of permutations on $n$ elements; the map extends linearly.
We give then the following definition: we denote as 
\begin{equation*}
\sym{\cdot,\cdot}:\si{k}\otimes\si{h}\xymatrix{\ar[r]&} \si{h+k}
\end{equation*}
the symmetrization of the two factors, more explicitly if
\begin{equation*}
\alpha=\sum_{\pi^{k}}\alpha_{\pi^{k}(1)}\otimes\cdots\otimes\alpha_{\pi^{k}(k)}\in\si{k}\quad,\quad\beta=
\sum_{\pi^{h}}\beta_{\pi^{h}(1)}\otimes\cdots\otimes\beta_{\pi^{h}(h)}\in\si{h}
\end{equation*}
then
\begin{equation*}
\sym{\alpha\otimes\beta}=\sum_{\pi^{k},\pi^{h}}S\big(\alpha_{\pi^{k}(1)}\otimes\cdots\otimes\alpha_{\pi^{k}(k)}
\otimes\beta_{\pi^{h}(1)}\otimes\cdots\otimes\beta_{\pi^{h}(h)}\big)\,. 
\end{equation*}
In particular we denote by $\sigma$ the map $\sym{\cdot,\cdot}$ when the first index is $1$:
\begin{equation}\label{defsigma}
\sigma:=\sym{\cdot,\cdot}:\xymatrix{\si{1}\otimes\si{i}\ar[r]&\si{i+1}}\,.
\end{equation}

Consider now for simplicity the case that $E$ corresponds to an irreducible $Sp(1)$ representation; then
\begin{equation*}
T_{x}M_{\C}\cong \si{1}\otimes\si{i-1}
\end{equation*}
and using Clebsch-Gordan relation,\,we obtain
\begin{equation}\label{s1s2immers}
T_{x}M_{\C}\cong \si{i}+\si{i-2}\xymatrix{\ar@{^{(}->}[r]&}\si{i+2}+\si{i}+\si{i-2}\cong \si{2}\otimes\si{i}\,;
\end{equation}
more precisely $T_{x}M_{\C}$ coincides with the kernel of the symmetrization 
\begin{equation*}
\sym{\cdot,\cdot}:\si{2}\otimes\si{i}\xymatrix{\ar[r]&} \si{i+2}\,.
\end{equation*}

{\bfseries Example.} There are (up to conjugation) three non-trivial homomorphisms $Sp(1)\to Sp(2)$: two correspond to the roots,\,
but in these cases the decomposition of the standard $Sp(2)$ representation $\C^{4}$ is not irreducible; in fact
\begin{equation*}
E=\C^{4}=\si{0}+\si{0}+\si{1}
\end{equation*}
for the long root, and comparing with the known decomposition of the adjoint representation one has
\begin{equation*}
\spd=S^{2}(\C^{4})=S^{2}(2\si{0}+\si{1})=\si{2}+2\si{1}+3\si{0}\,;
\end{equation*} 
for the short root we have instead
\begin{equation*}
E=\C^{4}=\si{1}+\si{1}
\end{equation*}
as in fact
\begin{equation*}
\spd=S^{2}(\C^{4})=S^{2}(2\si{1})=3\si{2}+\si{0}\,.
\end{equation*} 
There is a third embedding, corresponding to the $\mathfrak{sl}(2,\C)$ triple
\begin{align*}
 X&=\begin{pmatrix}
0&\sqrt{3}&0&0\\
0&0&0&\sqrt{2}\\
0&0&0&0\\
0&0&-\sqrt{3}&0 
\end{pmatrix}\,,
Y=\begin{pmatrix}
0&0&0&0\\
\sqrt{3}&0&0&0\\
0&0&0&-\sqrt{3}\\
0&\sqrt{2}&0&0 
\end{pmatrix}\,,\\
&\quad\quad\quad\quad\qquad H=\begin{pmatrix}
3&0&0&0\\
0&1&0&0\\
0&0&-3&0\\
0&0&0&-1
\end{pmatrix}\,,
\end{align*}
obtained using the recipe in \cite{collmcg},\,for which
\begin{equation}\label{irredspuspd}
E=\C^{4}=\si{3}\,.
\end{equation}

 \noindent\textit{Observation.} This last can be interpreted in the following way: recall that the decomposition of the 
Lie algebra $\la_{2}$ with respect to $\soq\subset\la_{2}$ is given by
\begin{equation*}
\si{2}_{+}+\si{2}_{-}+\si{1}_{-}\otimes\si{3}_{+}\,,
\end{equation*}
where $\si{k}_{\pm}$ denote the representations of the $\spu$ corresponding to the long ($+$) or to the short($-$)
root;\,so considering the diagonal embedding 
\begin{equation*}
\xymatrix{\spu_{\Delta}\ar@{^{(}->}[r]&\soq=\spu_{+}+\spu_{-}\ar@{^{(}->}[r]&\spu_{+}+\spd}\,,
\end{equation*}
consistently with the $Sp(1)Sp(2)$ structure of the Wolf space 
\begin{equation*}
\frac{G_{2}}{SO(4)}\,,
\end{equation*}
we have a description of its tangent space in the $EH$ formalism as $H\otimes E\cong\si{1}\otimes \si{3}$,\,corresponding to the 
representation in (\ref{irredspuspd}).
\vspace{3mm}

The action of $S^{2}H\cong \si{2}$ on $T_{x}M_{\C}$ can be therefore expressed suitably exploiting this new
formulation,\,involving the $\si{2}$ factor instead of the $\si{1}=H$; to understand more deeply this 
$\si{2}$-approach we need to define more explicitly the invariant immersion (\ref{s1s2immers}).\,
Let us define the map as
\begin{equation}\label{defQ}
Q:\xymatrix{\si{1}\otimes\si{i-1}\ar[r]^-{\omega_{H}\otimes}&\underline{\si{1}}\otimes\underline{\underline{\si{1}}}
\otimes\underline{\si{1}}\otimes\underline{\underline{\si{i-1}}}\ar[r]^-{\sym{\cdot,\cdot}}& \si{2}\otimes\si{i}}
\end{equation}
acting in the following way:\,if 
\begin{equation*}
Y=h\otimes \beta+\hh \otimes\hbeta\in \si{1}\otimes\si{i-1},\quad \beta,\hbeta\in \si{i-1}
\end{equation*}
then 
\begin{equation}\label{defQY}
Q(Y)=\frac{1}{2} \sym{h h}\sym{\hh\beta}+\frac{1}{4}(h\hh+\hh h)\big(\sym{\hh\hbeta}-\sym{h\beta}\big)-\frac{1}{2}\sym{\hh\hh}\sym{h\hbeta}
\end{equation}
is obtained, after tensorizing with the invariant element $\omega_{H}$, by symmetrization of the tensorial factors 
in accordance with the simple or double underlining marks in (\ref{defQ}).

Our next aim is to express the quaternionic action in terms of this description: a first guess in this sense
is that for $Q(Y)=\sum v_{i}\otimes p_{i}$ then
\begin{equation*}
Q(I_{1}Y)=v_{2}\otimes p_{3}+v_{3}\otimes p_{2}\,,
\end{equation*}
mimicking the adjoint representation of $\sud$ on itself;\,but this is not correct, as at the second step
\begin{equation*}
Q(I^{2}_{1}Y)=-v_{2}\otimes p_{2}-v_{3}\otimes p_{3}\,,
\end{equation*}
which is not $-Id$.\,Something more is needed to ``reconstruct" the missing term $-v_{1}\otimes p_{1}$.\,

The next Proposition gives the correct answer in order to express the quaternionic action from the $\si{2}$ viewpoint: 
\begin{prop}
Let $Y\in T_{x}M=\si{1}\otimes\si{i}$;\,if $Q(Y)=\sum v_{i}\otimes p_{i}$ then
\begin{equation}\label{quatalgpofview}
Q(I_{1}Y)=v_{1}\otimes \,\frac{1}{4}\sigma(Y)+v_{2}\otimes p_{3}-v_{3}\otimes p_{2}\,.
\end{equation}
\end{prop}
 \noindent\textit{Proof.} We have the definition of $Q(Y)$ as in (\ref{defQY}): then if we identify
$v_{i}$ with the basis $I_{i}$ defined in (\ref{quatbasis}),\,grouping the terms properly we obtain
\begin{eqnarray*}
p_{1}&=&-\frac{\imath}{4}\big(\sym{\hh\hbeta}-\sym{h\beta}\big)\\
p_{2}&=&\frac{1}{4}\big(\sym{\hh\beta}-\sym{h\hbeta}\big)\\
p_{3}&=&-\frac{\imath}{4}\big(\sym{\hh\beta}+\sym{h\hbeta}\big)   \,;
\end{eqnarray*}
the quaternionic action of $I_{1}$ on $Y$ is given,\,in the $\si{1}$ context,\,by
\begin{equation*}
I_{1}Y=-\imath h\otimes\beta+\imath \hh\otimes\hbeta\,;
\end{equation*}
so we obtain
\begin{equation*}
Q(I_{1}Y)=-\frac{\imath}{2}\sym{hh}\sym{\hh\beta}+\frac{\imath}{4}(h\hh+\hh h)\big(\sym{h\beta}+\sym{\hh\hbeta}\big)
-\frac{\imath}{2}\sym{\hh\hh}\sym{h\hbeta}
\end{equation*}
and in the form $Q(I_{1}Y)=\sum_{i=1}^{3}v_{i}\otimes q^{1}_{i}$ we have
\begin{eqnarray*}
q^{1}_{1}&=&\frac{\imath}{4}\big(\sym{h\beta}+\sym{\hh\hbeta}\big)\\
q^{1}_{2}&=&-\frac{i}{4}\big(\sym{\hh\beta}+\sym{h\hbeta}\big)\\
q^{1}_{3}&=&-\frac{1}{4}\big(\sym{\hh\beta}-\sym{h\hbeta}\big)  \,;
\end{eqnarray*}
the conclusion follows by the definition of $\sigma$ (\ref{defsigma}) and comparing the two sets of equalities.\,$\blacksquare$

In the same way we obtain for the other quaternionic elements
\begin{eqnarray*}
I_{2}Y&=&-\hh\otimes\beta+h\otimes\hbeta \\
I_{3}Y&=&\imath\hh\otimes\beta+\imath h\otimes\hbeta
\end{eqnarray*}
so that
\begin{eqnarray*}
Q(I_{2}Y)=&\frac{1}{2}\sym{hh}\sym{\hh\hbeta}-\frac{1}{2}\sym{h\hh}\big(\sym{\hh\beta}+\sym{h\hbeta}\big)
+\frac{1}{2}\sym{\hh\hh}\sym{h\beta}\\
Q(I_{3}Y)=&\frac{\imath}{2}\sym{hh}\sym{\hh\hbeta}+\frac{\imath}{2}\sym{h\hh}\big(\sym{\hh\beta}-\sym{h\hbeta}\big)
-\frac{\imath}{2}\sym{\hh\hh}\sym{h\beta}
\end{eqnarray*} 
which imply the equalities
\begin{equation*}
q^{i}_{j}=\eta_{ijk}\,p_{k}-\delta_{i}^{j}\,\frac{1}{4}\,\sigma(Y)\,,
\end{equation*}
where $\eta_{ijk}=sgn_{ijk}$ if $i\neq j$,\,otherwise $\eta_{iik}=0$;\,moreover
\begin{equation*}
p_{i}=-\frac{1}{4}\sigma(I_{i}Y)\,.
\end{equation*}
We can therefore state the quaternionic
relations in terms of this description: for example
\begin{align*}
Q(I_{1}^{2}Y)=Q(I_{1}I_{1}Y)&=-v_{1}\otimes\frac{1}{4}\sigma(I_{1}Y)-v_{2}\otimes p_{2}-v_{3}\otimes p_{3} \notag\\
&=-v_{1}\otimes p_{1} -v_{2}\otimes p_{2}-v_{3}\otimes p_{3}\notag\\
&=-Q(Y)
\end{align*}
and also
\begin{align*}
Q(I_{1}I_{2}Y)&=-v_{1}\otimes\frac{1}{4}\sigma(I_{_2}Y)-v_{2}\otimes q_{3}^{2}-v_{3}\otimes q_{2}^{2} \notag\\
&=-v_{1}\otimes p_{2} +v_{2}\otimes p_{1}-v_{3}\otimes \frac{1}{4}\sigma(Y)\notag\\
&=Q(I_{3}Y)
\end{align*}
as expected.



\section{The Coincidence Theorem}\label{thecoincidencetheorem}

Another way of expressing the \emph{twistor equation} (\ref{twisteqnqk}) is given by 
\begin{equation}\label{twisteq4}
\conn{S^{2}H}\mu_{A}=k\sum_{i=1}^{3}I_{i}\tilde{A}^{\flat}\otimes I_{i}\,,
\end{equation}
where $\tilde{A}$ is the Killing vector field generated by $A$ in $\la$, the symbol $\flat$ means passing to the corresponding
$1$-form via the metric and $k$ is the scalar curvature, which is constant as the metric is Einstein (for simplicity we can put $k=1$). 
On the other hand on $\taut$ we have defined the sections $s_{A}$ and the natural connection $\conn{\taut}$ so that (see (\ref{vperpingr}) 
and Proposition \ref{costantprojsec})
\begin{equation*}
\conn{\taut} s_{A}=\sum_{i=1}^{3} s^{\perp}_{A}\otimes v_{i}\otimes v_{i}\,.
\end{equation*}
\vspace{2mm}

In general,\,given a differentiable map $\Psi:M\to N$ of manifolds, and an isomorphism $\mlift$ between vector bundles 
$E\to F$ on the manifold $M$ and $N$ respectively,\,the second one equipped with a connection $\conn{F}$, 
we can define the \emph{pullback connection} $\hpsi^{*}\conn{F}$ acting in the following way on elements $\sect$ of $\Gamma(E)$:
\begin{equation*}
(\Psi^{*}\conn{F})_{Y}(\sect):=\hpsi^{*}(\conn{F}_{(\Psi_{*}Y)}(\hpsi \sect))
\end{equation*}
where $Y\in T_{x}M$ and $\hpsi^{*}$ means taking the pullback section.  

We want to apply this construction in our case, with the map $\Psi:M\to\graor(\la)$ induced by $\mu$, $N=\graor(\la),\,E=S^{2}H,\,F=\taut$;
our aim is to relate, at a fixed point $x\in M$, the action of the quaternionic structure on $1$-forms induced by $G$ (the duals of 
the Killing vector fields) with special cotangent vectors on the Grassmannian $\graor(\la)$:
\vspace{2mm}

\begin{lemma}\label{coincidencethm}
Let $M,\la,\graor(\la),\mu$ be defined as usual, with
\begin{equation*}
\mu=\sum_{i=1}^{3} I_{i}\otimes B_{i}
\end{equation*}
where $B_{i}=\conf v_{i}$,\,$\conf$ a differentaiable $G$-invariant function on $M$ and $v_{i}$ an orthonormal 
basis of a point $V\in\graor(\la)$;\,let us choose $A\in V^{\perp}\subset\la$;\,then at the point $x$ such that
$\Psi(x)=V$,\,for $\Psi$ induced by $\mu$ as usual,\,we have
\begin{equation}\label{coinceqn}
\frac{1}{\conf}\,I_{i}\tilde{A}^{\flat}=\Psi^{*}(A\otimes v_{i})^{\flat}\,,
\end{equation}
where $A\otimes v_{i}\in T_{x}\graor(\la)$.\,Moreover we have $\|\mu\|^{2}=3\conf^{2}$.
\end{lemma}

 \noindent\textit{Proof.} Let $\Psi$ denote the conformal lift of the map $\mu$ so that 
\begin{equation}
\Psi(I_{i})=\frac{1}{\conf^{2}}B_{i}\,;
\end{equation}
hence as seen in Proposition \ref{proptwisteq}
\begin{equation*}
\Psi(\mu_{A})=s_{A}\,;
\end{equation*}
then applying the $\Psi^{*}\conn{\taut}$ connection of $S^{2}H$ to $\mu_{A}$ we obtain
\begin{align}\label{connpullmu}
(\Psi^{*}\conn{\taut})\mu_{A}&=\Psi^{*}\left(\conn{\taut}(\Psi (\mu_{A}))\right) \notag\\
&=\Psi^{*}\left(\conn{\taut}s_{A}\right)\notag\\
&=\Psi^{*}\left(\sum_{i=1}^{3}s_{A}^{\perp}\otimes v_{i}\otimes v_{i}\right)\notag\\
&=\conf \sum_{i=1}^{3}\Psi^{*}(s_{A}^{\perp}\otimes v_{i})\otimes I_{i}\,;
\end{align}
on the other hand the difference of two connections on the same vector bundle is a tensor, so given any section 
$\sect\in S^{2}H$ which vanishes at a point $x\in M$ 
\begin{equation*}
(\conn{S^{2}H}-\Psi^{*}\conn{\taut})\sect (x)=0\,.
\end{equation*}
This is precisely the case for the section $\mu_{A}$ at the point $x$ for which $\Psi(S^{2}H_{x})=V$, because
$A\in V^{\perp}$ by hypothesis; in other words
\begin{equation*}
\conn{S^{2}H}\mu_{A}\,_{|_{x}}=(\Psi^{*}\conn{\taut})\mu_{A}\,_{|_{x}}\,.
\end{equation*}
In the light of the calculations in (\ref{connpullmu}) and of the twistor equation (\ref{twisteq4}), we can deduce
\begin{equation*}
\sum_{i=1}^{3}I_{i}\tilde{A}^{\flat}\otimes I_{i}=\conf \sum_{i=1}^{3}\Psi^{*}(s_{A}^{\perp}\otimes v_{i})\otimes I_{i}\,;
\end{equation*}
the result follows considering that $s_{A}^{\perp}=A$ at $V$. $\blacksquare$\\
\vspace{2mm}

Lemma \ref{coincidencethm} leads to various ways of relating elements in the respective spaces $T_{x}M$ and
$T_{V}\graor(\la)$ and the quaternionic elements $I_{i}$; nevertheless it is stated merely in terms of $1$-forms, whereas we
are interested in involving the two metrics in this interplay. To this aim, let us define a linear transformation
$\natural$ of $T_xM$ by \begin{equation}\label{naturaldef} \autnatural{X}:=\big(\Psi^{*}\big((
\Psi_{*}X)^{\flat}\big)\big)^{\sharp} \end{equation} in $\vend(T_{x}M)$. This corresponds to moving in a counterclockwise
sense around the following diagram, starting from bottom left: \begin{equation}\label{noncommdiag}
\xymatrix{T^{*}_{x}M\ar[d]_-{\sharp}&T^{*}_{V}\gra\ar[l]_-{\Psi^{*}}\\ T_{x}M\ar[r]_-{\Psi_{*}}& T_{V}\gra\ar[u]_-{\flat}
}\,.  \end{equation} Thus the linear endomorphism $\autnatural{\cdot}$ measures the noncommutativity of the diagram
(\ref{noncommdiag}), and the difference between the pullbacked Grassmannian metric from the quaternionic one.
\vspace{2mm}

We are in position now to prove the \emph{Coincidence Theorem}:

\begin{teor}\label{coinccor}
Let $Y\in T_{x}M$ such that 
\begin{equation*}
\Psi_{*}Y=\sum v_{i}\otimes p_{i}\,;
\end{equation*}
for $p_{i}\in V^{\perp}$ with $V=\Psi(x)$; then 
\begin{equation*}
\autnatural{Y}=\frac{1}{\conf}\sum_{i}I_{i}\tilde {p}_{i}.
\end{equation*}
\end{teor}
 \noindent\textit{Proof.} Using the definitions and (\ref{coinceqn}) we obtain
\begin{align*}
(\Psi_{*}Y)^{\flat}\big(\Psi_{*}Z\big)&=\langle \sum v_{i}\otimes p_{i}\,,\, \Psi_{*}Z\rangle_{\graor}\notag\\
&=\frac{1}{\conf}\langle \sum I_{i}\tilde{p}_{i}\,,\, Z\rangle_{M}
\end{align*}
for any $Z\in T_{x}M$,\,hence the conclusion.\,$\blacksquare$
\vspace{2mm}

 Theorem \ref{coinccor} provides a memorable way of ``converting" tangent vectors of $\gra(\la)$ to tangent vectors on $M$ by means
of the correspondence \begin{align*} &\xymatrix{v_{i}\ar[r]&I_{i}}\notag\\ &\xymatrix{p_{i}\ar[r]&\tilde{p_{i}}}\notag
\end{align*} for $p_{i}\in V^{\perp}$. 
\vspace{2mm}

The equivariance of the moment map $\mu$ implies that Killing vector fields on $M$ map to Killing vector fields on
$\graor(\la)$: in other words if $\tilde{A}$ is induced by $A\in\la$ on $M$,\,then \begin{equation*}
\Psi_{*}\tilde{A}=\sum_{i=1}^{3}v_{i}\otimes [A\,,\,v_{i}]^{\perp}\,.  \end{equation*} Let now
$\alpha=(\sum_{i=1}^{3}v_{i}\otimes p_{i})^{\flat}\in T_{x}^{*}\graor(\la)$ and let $A_{r}$ be an orthonormal basis of
$V^{\perp}$;\,then \begin{align*} \sum _{r=1}^{n-3}\langle \Psi^{*}\alpha,\,\tilde{A}_{r}\rangle A_{r}&=\sum
_{r=1}^{n-3}\langle \alpha,\,\Psi_{*}\tilde{A}_{r}\rangle A_{r} =\sum _{i,r}\langle p_{i},\,[v_{i}\,,\,A_{r}]^{\perp}\rangle
A_{r}\notag\\ &=\sum _{i,r}\langle \,p_{i},\,[v_{i}\,,\,A_{r}]\,\rangle A_{r}=\sum _{i,r}\langle
[p_{i},\,v_{i}]\,,\,A_{r}\rangle A_{r}\notag\\ &=\sum _{i} [p_{i},\,v_{i}]^{\perp}\,.  \end{align*} We can therefore define a
mapping \begin{equation}\label{rhodef} \rho:\xymatrix{T_{x}^{*}M\ar[r] & V^{\perp}_{_{}}} \end{equation} by
$\rho(\ufd)=\sum_{r}\langle \ufd\,,\, \tilde{A}_{r}\rangle A_{r}$; so if $\alpha\in T_{x}^{*}\graor(\la)$,\,then
$\Psi^{*}\alpha \in T_{x}^{*}M$,\,and the composition $\tilde{\gop}=\rho\circ \Psi^{*}$ is a map \begin{equation*}
\tilde{\gop}:\xymatrix{T_{x}^{*}\graor(\la) \ar[r] & V^{\perp}}\,.  \end{equation*} defined by
$\tilde{\gop}(\alpha)=\sum_{i}[v_{i}\,,\,p_{i}]^{\perp}$; this operator can be described as \begin{equation*}
\tilde{\gop}=\pi^{\perp}\circ\gop\, \end{equation*} where $\gop(\alpha)=\sum_{i}[v_{i}\,,\,p_{i}]$ is the obstruction to the orthogonality of $\alpha$ 
to the $G$-orbit: in fact

\begin{lemma}\label{orthogobstr}
A tangent vector $P=\sum_{i=1}^{3}v_{i}\otimes p_{i}\in T_{V}\graor(\la)$ is orthogonal to the $G$-orbit through the 
point $V$ if and only if $\gop(P)=0$.
\end{lemma}
 \noindent\textit{Proof.} For any $A\in \la$ let us consider the Killing vector field $\tilde{A}$ on $\graor(\la)$; the condition 
of orthogonality of $P$ is expressed by
\begin{align*}
0&=\langle\, \tilde{A}\,,\, P\,\rangle=\sum_{i=1}^{3}\langle\,[A\,,\,v_{i}]^{\perp}\,,\,p_{i}\,\rangle=\notag\\
&=\sum_{i=1}^{3}\langle\,[A\,,\,v_{i}]\,,\,p_{i}\,\rangle
=\sum_{i=1}^{3}\langle\, A\,,[\,v_{i}\,,\,p_{i}]\,\rangle=\notag\\
&=\langle\, A\,,\, \gop(P)\,\rangle\quad.\quad\blacksquare
\end{align*}
\vspace{2mm}

We give now a more explicit description of the quaternionic endomoprhisms:

\begin{prop}\label{coinctmcons}
Let $Y\in T_{x}M$ so that
\begin{equation*}
\Psi_{*}Y=v_{1}\otimes p_{1}+ v_{2}\otimes p_{2}+ v_{3 }\otimes p_{3}\,;
\end{equation*}
then we have
\begin{equation}\label{quatexprct}
\Psi_{*}I_{1}Y=\frac{1}{\conf}\, v_{1}\otimes \rho(Y^{\flat})-\,v_{2}\otimes p_{3}+ v_{3}\otimes p_{2}\,.
\end{equation}
\end{prop}
 \noindent\textit{Proof.} Consider any $A\in V^{\perp}$,\,then
\begin{align}\label{p1proj}
\langle p_{1}\,,\, A\rangle _{K}&=\langle \Psi_{*}Y\,,\, A\otimes v_{1}\rangle_{\gra}=\frac{1}{\conf}\langle I_{1}\tilde{A}^{\flat}\,,\, Y \rangle\notag\\
&=\frac{1}{\conf}\langle I_{1}\tilde{A}\,,\, Y \rangle_{M}=-\frac{1}{\conf}\langle \tilde{A}\,,\, I_{1} Y \rangle_{M}\notag\\
&=-\frac{1}{\conf}\langle I_{1} Y^{\flat} \,,\, \tilde{A} \rangle\,,
\end{align}
where $\langle\,,\,\rangle_{M,\mathbb{G}}$ denote the respective Riemannian metrics, $\langle\,,\,\rangle_{K}$ minus the Killing 
form on $\la$ and $\langle\,,\,\rangle$ without subscript
is merely the contraction of a cotangent and tangent vector; then considering (\ref{p1proj}) and (\ref{rhodef})
\begin{align*}
p_{1}&=\sum_{r}\langle p_{1}\,,\, A_{r}\rangle _{K}\,A_{r}=-\frac{1}{\conf}\sum_{r}\langle I_{1} Y^{\flat} \,,\, \tilde{A}_{r} \rangle\, A_{r}\notag\\
&=-\frac{1}{\conf}\rho(I_{1}Y^{\flat})
\end{align*}
and analogously
\begin{equation*}
p_{i}=-\frac{1}{\conf}\rho(I_{i}Y^{\flat})\,,\quad i=2,3\,;
\end{equation*}
in consequence
\begin{align*}
\Psi_{*}I_{1}Y&=\frac{1}{\conf}\, v_{1}\otimes \rho(Y^{\flat})-\frac{1}{\conf} \,v_{2}\otimes \rho(I_{3}Y^{\flat})+
\frac{1}{\conf} \,v_{3}\otimes \rho(I_{2}Y^{\flat})\notag\\
&=\frac{1}{\conf}\, v_{1}\otimes \rho(Y^{\flat})-\,v_{2}\otimes p_{3}+ v_{3}\otimes p_{2}\,.\,\blacksquare
\end{align*}
\vspace{3mm}

Clearly analogous assertions are valid for $I_{2}$ and $I_{3}$. 
\vspace{2mm}

 \noindent\textit{Remark.} Assuming that $\Psi_{*}$ is injective at the point $x$, we can define the 
\emph{push forward} of the endomorphisms $I_{k}$ in the obvious way, namely via the equation:
\begin{equation*}
(\Psi_{*}I_{k})Z:= \Psi_{*}(I_{k}(\Psi_{*}^{-1}Z)
\end{equation*} 
and Proposition \ref{coinctmcons}. A striking feature of (\ref{quatexprct}) is that in the expression obtained the first summand is independent from 
$I_{1}$. The operators $\rho,\,\gop$ appear as the essential ingredient to reconstruct the quaternionic action; the other
summands $-v_{2}\otimes p_{3}+v_{3}\otimes p_{2}$ are obtained from the adjoint representation and (as explained in Section \ref{thespnsp1structure}) 
are not sufficient. 
Nevertheless proposition \ref{coinctmcons} predicts that if $Y$ is perpendicular to the $G$-orbit on 
$M$,\,then 
\begin{equation*}
\rho(Y^{\flat})=0\,,
\end{equation*}
thanks to the definition of $\rho$ (see Lemma \ref{orthogobstr}); in that case
\begin{equation*}
\Psi_{*}I_{1}Y=-\,v_{2}\otimes p_{3}+ v_{3}\otimes p_{2}
\end{equation*}
which coincides with the irreducuble representation of $\spu$ on $V=\R^{3}$.


\section{Examples and applications}\label{examplesandapplications}

 The apparent distinction between the points of view we have adopted in Section \ref{thespnsp1structure} and Section
\ref{thecoincidencetheorem} disappears as soon as one compares (\ref{quatexprct}) and (\ref{quatalgpofview}). This suggests
that an intimate relationship exists between the two descriptions of the quaternionic structure: we are going to discuss now
some examples which throw light on this link.

Let us consider the Wolf space 
\begin{equation*}
\Ha\Pro^{1}\cong \frac{Sp(2)}{Sp(1)\times Sp(1)}\cong \frac{SO(5)}{SO(4)}\cong S^{4}
\end{equation*}
and the action of the stabilizer $Sp(1)\times Sp(1)$ of a point $N$, with Lie algebra $\spu_{+}\oplus \spu_{-}=\soq$; 
this is a cohomogeneity $1$ action, with generic orbits isomorphic to 
\begin{equation*}
S^{3}\cong \frac{Sp(1)\times Sp(1)}{Sp(1)_{\Delta}}
\end{equation*}
where $Sp(1)_{\Delta}$ is the diagonal representation, and $2$ singular orbits
corresponding to a couple of antipodal points $N, S$. Let us choose at the point $N$ any closed geodesic
$\beta(t)$ connecting $N$ to $S$: this will be orthogonal to any $Sp(1)\times Sp(1)$ orbit, and will intersect 
all of them (a \emph{normal geodesic} in the language of \cite{alek-alek}, which in higher cohomogeneity 
is generalized by submanifolds called \emph{sections}, see \cite{hptt1}). For instance, we can 
choose $N= e\,Sp(1)\times Sp(1)$, and take the geodesic corresponding to following copy of $U(1)\subset Sp(2)$:
\begin{equation}\label{tanggeod}
g(t)=
\begin{pmatrix}
\cos{t} & \sin{t} & 0 & 0\\
-\sin{t} & \cos{t} & 0 & 0 \\
0 & 0 & \cos{t} & \sin{t}\\
0 & 0 & -\sin{t} & \cos{t} \\
\end{pmatrix}=\exp
\begin{pmatrix}
0 & t & 0 & 0\\
-t & 0 & 0 & 0 \\
0 & 0 & 0 & t\\
0 & 0 & -t & 0 \\
\end{pmatrix}
\,,
\end{equation}
where the matrix on the right is denoted by $t\, u$. This subgroup generates 
a geodesic $\beta(t)$ connecting $N\; (t=0)$ with the south pole $S\; (t=\pi/2)$ passing through the equator $(t=\pi/4)$,\,and 
then backwards to $N\;(t=\pi)$. The stabilizer of the $Sp(1)\times Sp(1)$ action is constant along $\beta(t)$ on points that are 
different from  $N$ and $S$, and coincides with $Sp(1)_{\Delta}$, both along $\beta(t)$ in $\Ha\Pro^{1}$ and along $\mathfrak{u}(1)$ for the isotropy
representation. 

Let now $e_{i}$ and $f_{i}$ denote orthonormal bases of $\spu_{+}$ and $\spu_{-}$ respectively; as $\soq$ is a subalgebra of
$\spd$ corresponding to the longest root,\,the elements of the two copies of $\spu$ correspond to the following
matrices:
\begin{equation}\label{basis1}
e_{1}=\frac{1}{\sqrt{2}}\begin{pmatrix}
\imath & 0 & 0 & 0\\
0 & 0 & 0 & 0 \\
0 & 0 & -\imath & 0\\
0 & 0 & 0 & 0 \\
\end{pmatrix}\quad,\quad
f_{1}=\frac{1}{\sqrt{2}}\begin{pmatrix}
0 & 0 & 0 & 0\\
0 & \imath & 0 & 0 \\
0 & 0 & 0 & 0\\
0 & 0 & 0 & -\imath \\
\end{pmatrix}\,,
\end{equation}

\begin{equation}\label{basis2}
e_{2}=\frac{1}{\sqrt{2}}\begin{pmatrix}
 0 & 0 & 1 & 0\\
0 & 0 & 0 & 0 \\
-1 & 0 & 0 & 0\\
0 & 0 & 0 & 0 \\
\end{pmatrix}\quad,\quad
f_{2}=\frac{1}{\sqrt{2}}\begin{pmatrix}
0 & 0 & 0 & 0\\
0 & 0 & 0 & 1 \\
0 & 0 & 0 & 0\\
0 & -1 & 0 & 0 \\
\end{pmatrix}\,
\end{equation}

\noindent and

\begin{equation}\label{basis3}
e_{3}=\frac{1}{\sqrt{2}}\begin{pmatrix}
0 & 0 & \imath & 0\\
0 & 0 & 0 & 0 \\
\imath & 0 & 0 & 0\\
0 & 0 & 0 & 0 \\
\end{pmatrix}\quad,\quad
f_{3}=\frac{1}{\sqrt{2}}\begin{pmatrix}
0 & 0 & 0 & 0\\
0 & 0 & 0 & \imath \\
0 & 0 & 0 & 0\\
0 & \imath & 0 & 0 \\
\end{pmatrix}\,;
\end{equation}
so if $e_{i}(t)$ and $f_{i}(t)$ denote an orthonormal basis of the isotropy subalgebra at $\beta(t)$ (given
by $Ad_{g(t)}\soq$), we get via the Killing metric:
\begin{align*}
\langle e_{i}\,,\,f_{j}(t)\rangle&=\delta^{i}_{j}\sin^{2}{t}\\
\langle e_{i}\,,\,e_{j}(t)\rangle&=\delta^{i}_{j}\cos^{2}{t}\\
\langle f_{i}\,,\,e_{j}(t)\rangle&=\delta^{i}_{j}\sin^{2}{t}\\
\langle f_{i}\,,\,f_{j}(t)\rangle&=\delta^{i}_{j}\cos^{2}{t}\,;
\end{align*}
in terms of Killing vector fields this implies
\begin{equation*}
\pi_{S^{2}H}(\nabla\tilde{e_{i}})=\sin^{2}{t}\;f_{i}(t)\quad,\quad\pi_{S^{2}H}(\nabla\tilde{f_{i}})=\cos^{2}{t}\;f_{i}(t)\,.
\end{equation*}
if we identify $S^{2}H\cong Ad_{g(t)}\spu_{-}$. 

The conclusion is that along $\beta(t)$ the moment map for the action of $Sp(1)\times Sp(1)$ on $\Ha\Pro^{1}$ 
is given by
\begin{equation}\label{mmapsoq}
\mu(\beta(t))=\sum_{i}\omega_{i}\otimes\,(\cos^{2}{t}\;f_{i}+\sin^{2}{t}\;e_{i})\,,
\end{equation}
up to a constant. This is the only information that we need to reconstruct the moment map on the whole $\Ha\Pro^{1}$, 
as $\beta(t)$ intersects all the orbits and the moment map is equivariant.

We can now interpret these facts in terms of the induced map  
\begin{equation*}
\Psi:\Ha\Pro^{1}\xymatrix{\ar[r]&}\graor (\soq)\,;
\end{equation*}
first of all we note that in this case $M_{0}=M$, 
as the three vectors 
\begin{equation}\label{psispan}
B_{i}(t)= \cos^{2}{t}\;f_{i}+\sin^{2}{t}\;e_{i}
\end{equation}
are linearly independent for all $t$; moreover we observe that $\hat{\Phi}$ is a conformal mapping of bundles, as asked
in the general hypotheses discussed in Section \ref{thetwotwistorequations}.

Recall from \cite{swann98} that the critical manifolds for the gradient flow of the functional 
\begin{equation*}
\psi=\langle[v_{1}, v_{2}], v_{3}\rangle
\end{equation*}
defined on $\graor (\soq)$ are given by the maximal points $\spu_{+},\,\spu_{-}$ and the submanifold
\begin{equation*}
C_{\Delta}=\R\Pro^{3}\cong \frac{Sp(1)\times Sp(1)}{\Z_{2}\times Sp(1)_{\Delta}}
\end{equation*}
corresponding to the $3$-dimensional subalgebra $\spu_{\Delta}$, for $\psi>0$; the unstable manifold $M_{\Delta}$ 
emanating from this last one is $4$-dimensional and isomorphic to 
\begin{equation*}
\frac{\Ha\Pro^{1}\setminus \{N,S\}}{\Z_{2}}\,.
\end{equation*}
A trajectory for the flow of $\nabla \psi$ is given by
\begin{equation}\label{trajpsi}
V(x,y)=\vspan\{x e_{i}+y f_{i} \; | \,x^{2}+y^{2}=1\,,\,i=1\cdots 3\}\,,
\end{equation}
therefore, comparing (\ref{trajpsi}) with (\ref{psispan}) we obtain that $\Psi(\Ha\Pro^{1})=M_{\Delta}\cup \spu_{+}\cup \spu_{-}$; in particular:
\begin{align}
\Psi(N)&=\spu_{-} \label{psi1} \\
\Psi(S)&=\spu_{+} \label{psi2} \\
\Psi(\beta(\pi/4))&=\spu_{\Delta}\,. \label{psi3}
\end{align}
\vspace{1mm}

 \noindent\textit{Observation.} The map $\Psi$ is not injective.\,The points corresponding to $t$ and $\pi-t$ are sent
to the same $3$-plane; so the principal orbits of type $S^{3}$ in $\Ha\Pro^{1}$ are sent to the orbits of type $\R\Pro^{3}$ 
in $M_{\Delta}$.\,The map $\Psi$ becomes injective on the orbifold $\Ha\Pro^{1}/\Z_{2}$, nevertheless $\Phi_{*}$ is
injective away from $N,\,S$.
\vspace{2mm}

Therefore the $Sp(1)\times Sp(1)$ orbit through $x_{\Delta}=\beta(\pi/4)$ is sent through $\Psi$ to the critical orbit $C_{\Delta}$; we
have

\begin{prop}
The differential
\begin{equation*}
T_{x_{\Delta}}\Ha\Pro^{1}\xymatrix{\ar[r]^-{\Psi_{*}}&}T_{\spu_{\Delta}}\graor(\soq)
\end{equation*}
is a linear $Sp(1)_{\Delta}$-invariant injective map. It coincides (up to a constant) with the map $Q$ defined in (\ref{defQ}), 
in terms of $Sp(1)_{\Delta}$ modules.
\end{prop}
 \noindent\textit{Proof.} Let $\alpha(t)$ be any curve through $x_{\Delta}$, then
\begin{equation*}
g_{*}\cdot\Psi_{*}\alpha'(0)=\frac{d}{dt}\,g\cdot\Psi(\alpha(t))=\frac{d}{dt}\,\Psi(g\cdot\alpha(t))=\Psi_{*}g_{*}\cdot\alpha'(0)
\end{equation*} 
for $g\in Sp(1)_{\Delta}\subset Sp(1)\times Sp(1)_{x_{\Delta}}$ where this last is the isotropy subgroup at $x_{\Delta}$; in this  
case the Lie algebra $\spu_{\Delta}$ of $Sp(1)_{\Delta}$, which is the stabilizer at $\beta(t)$ for any $t$, 
turns out to coincide with the image $\Psi(x_{\Delta})$. The decomposition of the holonomy representation in terms
of $Sp(1)_{\Delta}$-modules is given in this case by
\begin{equation*}
E\otimes H\cong \si{1}\otimes \si{1}\cong \si{2} +\si{0}\,;
\end{equation*} 
correspondingly, the decomposition of the Grassmannian's tangent space at $V=\spu_{\Delta}$ is given by
\begin{align*}
T_{V}\graor(\soq)&\cong V\otimes V^{\perp}\cong \spu_{\Delta}\otimes \si{2}\cong \si{2}\otimes\si{2}\notag\\
&\cong \si{4}+\si{2}+\si{0}\,,
\end{align*}
and $\Psi_{*}$ sends injectively $\si{2}+\si{0}$ in $ \si{4}+\si{2}+\si{0}$; then, as a consequence of Schur's lemma, 
an isomorphism of $Sp(1)_{\Delta}$-modules is unique up to a constant for each irreducible submodule, hence
\begin{equation*}
\begin{cases}
\Psi_{*}&=a \,Q \quad \text{on $\si{2}$}\\
\Psi_{*}&=b \,Q \quad \text{on $\si{0}$}
\end{cases}
\end{equation*}
for some constants $a,\,b\in \R$.\,$\blacksquare$
\vspace{2mm}

An analogous situation holds for appropriate orbits in the following cases, which are all cohomogeneity $1$ actions on 
classical Wolf spaces:
\begin{itemize}
\item $Sp(n)Sp(1)$ acting on $\Ha\Pro^{n}$;
\item $Sp(n)$ acting on $\grasc{2n}$;
\item $SO(n-1)$ acting on $\grasr{n}$.
\end{itemize}
In the first case the orbit sent through $\Psi$ to a critical submanifold of type $C_{\Delta}$ in the corresponding 
Grassmannian is one of the principal orbits $S^{4n-1}$, in the second and third case is one of the singular orbits, 
more precisely 
\begin{equation*}
\frac{Sp(n)}{Sp(n-2)\times U(2)}\quad\text{and}\quad\graor(\R^{n-1})\cong\frac{SO(n-1)}{SO(n-4)\times SO(3)}\,
\end{equation*} 
respectively.
\vspace{2mm}

This situation can be generalized in the following sense: let $G$ be a compact group acting by quaternionic
isometries on a QK manifold $M$; let $G_{x}$ denote the stabilizer at the point $x\in M$; 
then $G_{x}\subset SO(T_{x}M)$ with respect to the quaternionic metric. Since QK manifolds are carachterized by the
condition $\hol(M)_{x}\subset Sp(n)Sp(1)\subset SO(T_{x}M)$, we have by hypothesis that $G_{x}\subset Sp(n)Sp(1)$.
Now suppose that $G_{x}$ contains some copy of $Sp(1)$ with nontrivial ptojection on
the $Sp(n)$ factor. In the case that
\begin{equation*}
\Psi(x)=\spu
\end{equation*}
and that a tubular neighborhood of $G_{x}$ is sent to the unstable manifold (for $\psi>0$) emanating from the critical manifold
$C\subset \graor(\la)$ corresponding to $\spu$, then
we have a corresponding decomposition of $T_{x}M$ and $T_{\spu}\graor(\la)$ in $Sp(1)$-modules, and the differential
$\Psi_{*}$ coincides with $Q$ up to determining $2k$ constants, $2$ for each $Sp(1)$-irreducible summand of the standard 
$Sp(n)$ module $E$.

Let us now decompose the holonomy representation in the case that $Sp(1)$ is the standard quaternionic subgroup, 
hence with trivial projection on the $Sp(n)$ factor. In this case $E$ turns out to be a direct sum of trivial 
representations:
\begin{equation*}
E\otimes H \cong (2n\,\si{0})\otimes \si{1}\cong 2n\, \si{1}
\end{equation*}
where $2\si{1}$ can be identified with the complexified algebra of Quaternions $\Ha\otimes _{\R}\C$. Therefore
going back to the real tangent bundle, we obtain the $Sp(1)$ invariant decomposition
\begin{equation}\label{quatdecomp}
T_{x}M\cong n\, \Ha\,.
\end{equation}
The presence of the $G$-action allows to single out a quaternionic line of $T_{x}M$: this determines a quaternionic 
$1$-dimensional distribution $\nqd$ on $M$, or a section $\tau:M\xymatrix{\ar[r]&}\Ha\Pro(TM)$ of the associated $\Ha\Pro^{n-1}$-bundle.
\vspace{2mm}

The distribution $\nqd$ arises in the following way: recall that at a point $V\in \graor(\la)$ with $v_{1},\,v_{2},\,v_{3}$ 
ON basis, we have 
\begin{equation*}
\grad\,\psi= v_{1}\otimes [v_{2},v_{3}]^{\perp}+v_{2}\otimes [v_{3},v_{1}]^{\perp}+v_{3}\otimes [v_{1},v_{2}]^{\perp}.
\end{equation*}
Maintaining the general hypotheses considered in Sections \ref{thetwotwistorequations} and \ref{thecoincidencetheorem}, and
assuming that $\Psi_{*}$ is injective, let us define $X:=\Psi_{*}^{-1}(\grad\,\psi)$; then we have:
\begin{lemma}
Suppose that $\Psi(x)=V$. Then the subspaces
\begin{eqnarray*}
\vspan\{\grad\,\psi,\,\tilde{v_{1}},\,\tilde{v_{2}},\,\tilde{v_{3}} \}&\subset& T_{V}\graor (\la)\\
\vspan\{X,\,\tilde{v_{1}},\,\tilde{v_{2}},\,\tilde{v_{3}} \}&\subset& T_{x}M
\end{eqnarray*}
are $Sp(1)$ invariant, hence quaternionic.
\end{lemma}
 \noindent\textit{Proof.} We need to prove that the endomorphisms of $S^{2}H$ over $x$ (or equivalently those of $\taut$ over $V$) 
preserve the respective subspaces; let us recall the description of $I_{1},\,I_{2},\,I_{3}$ given in Proposition 
\ref{coinctmcons}, then
\begin{align}\label{gradfv1}
I_{1}(\grad\,\psi)&=\frac{1}{\conf}\, v_{1}\otimes \rho\big((\grad\,\psi)^{\flat}\big)-\,v_{2}\otimes [v_{1},v_{2}]^{\perp}+ v_{3}\otimes [v_{3},v_{1}]^{\perp}\notag\\
&=-\,v_{2}\otimes [v_{1},v_{2}]^{\perp}+ v_{3}\otimes [v_{3},v_{1}]^{\perp}\notag\\
&=-\tilde{v_{1}}\,,
\end{align} 
where the first summand vanishes thanks to the $G$-invariance of $\psi$, which implies that $\grad\,\psi$ is orthogonal
to the $G$ orbits. Analogously, $I_{2}(\grad\,\psi)=-\tilde{v_{2}}$ and $I_{3}(\grad\,\psi)=-\tilde{v_{3}}$, and the 
quaternionic identities imply that the whole $\vspan\{\grad\,\psi,\,\tilde{v_{1}},\,\tilde{v_{2}},\,\tilde{v_{3}} \}$
is preserved; the second inclusion follows from the injectivity and equivariance of $\Psi$.\,$\blacksquare$
\vspace{2mm}

In all the examples discussed above the distribution $\nqd$ turns out to be integrable, with integral manifolds isomorphic
to $\Ha\Pro^{1}$ embedded quaternionically in $\Ha\Pro^{n}$, $\grasc{2n}$ or $\grasr{n}$ respectively. 

For $Sp(1)\times Sp(1)$ acting on $\Ha\Pro^{1}$ the ditribution $\nqd$ clearly coincides with the tangent bundle; in this
case it is possible to describe the relationship between the two metrics and the $\autnatural{\cdot}$ endomorphism:

\begin{prop}
Let $M=\Ha\Pro^{1}\setminus\{N,S\}$; consider the decomposition 
\begin{align}\label{hp1split}
T_{x}M&\cong \vspan\{\tilde{v_{1}},\tilde{v_{2}},\tilde{v_{3}}\}\oplus \vspan\{X\}\notag\\
&=:\summ_{1}\oplus \summ_{2}
\end{align}
induced by the $Sp(1)\times Sp(1)$ action; then the map $\Psi:M\xymatrix{\ar[r]&}\gra(\soq)$ satisfies the condition
\begin{equation}\label{confmap}
\Psi^{*}\langle\,,\,\rangle_{\gra}|_{\summ_{i}}=\eta_{i}(x)\langle\,,\,\rangle_{M}\,\quad i=1,\,2
\end{equation}
where $\eta_{i}(x)$ two real-valued $Sp(1) \times Sp(1)$ invariant function defined on $M$. The endomorphism 
(\ref{naturaldef}) is just the multiplication by $\eta_{i}(x)$ on $\summ_{i}$. 
\end{prop}

 \noindent\textit{Proof.} The tangent space $T_{V}\gra(\soq)$ along the unstable manifold can be seen as an irreducible $Sp(1)_{\Delta}$-module, and 
$\Psi_{*}$ as a morphism of $Sp(1)$-modules. Schur's Lemma guarantees the uniqueness of an invariant bilinear form 
(up to a constant), for every irreducible submodule. Recall that
\begin{equation*}
T_{x}M\cong \si{2}\oplus \si{0}
\end{equation*}
as $Sp(1)_{\Delta}$ representations, corresponding to the splitting (\ref{hp1split}): therefore equation (\ref{confmap}) 
holds, as both metrics are $Sp(1)_{\Delta}$ invariant. For the second assertion, let $Y\in \summ_{i}$:
\begin{align*}
\autnatural{Y}&=\big(\Psi^{*}\big(( \Psi_{*}Y)^{\flat}\big)\big)^{\sharp}\notag\\
&=\big(\Psi^{*}\big( \langle\Psi_{*}Y,\,\cdot\,\rangle_{\gra}\big)\big)^{\sharp}\notag\\
&=\eta_{i}(x)\big( \langle Y,\,\cdot\,\rangle_{M}\big)^{\sharp}\notag\\
&=\eta_{i}(x)Y\,
\end{align*}
as required. $\blacksquare$
\vspace{2mm}

 \noindent\textit{Observation.} Equation (\ref{gradfv1}) together with the equality $\|\grad\,\psi\|=3/2\,\|\tilde{v_{i}}\|$ confirms 
that the endomorphisms $I_{i}$ are \emph{not} isometries for Grassmannian metric; hence $\Psi^{*}\langle\,,\,\rangle_{\gra}$ 
and $\langle\,,\,\rangle_{M}$ can not coincide. Indeed,
\begin{equation*}
\|\grad\;\psi\|^{2}_{\gra}=\frac{3}{2}\,\|\tilde{v_{1}}\|^{2}_{\gra}=\frac{3}{2}\,\eta_{2}\,\|\tilde{v_{1}}\|^{2}_{M}\,;
\end{equation*}
moreover
\begin{equation*}
\|\grad\,\psi\|^{2}_{\gra}=\eta_{1}\,\|X\|^{2}_{M}
\end{equation*}
and $\|X\|_{M}=\|I_{1}X\|_{M}=\|\tilde{v_{1}}\|_{M}$. Thus $\frac{\eta_{1}}{\eta_{2}}=\frac{3}{2}$. An analogous result is 
expected to hold in general. 
\vspace{5mm}

{\footnotesize{\bfseries Acknowledgements.} This article is based on part of the author's PhD thesis at \emph{La Sapienza
}University of Rome, written under the supervision of S. Salamon, whom the author wishes especially to thank. He is also
grateful to A. F. Swann for useful discussions and comments.  }

\vspace{5mm}

\noindent\textsc{Dipartimento di Matematica ``G. Castelnuovo", Universit\ac ``La Sapienza", Piazzale A. Moro, 2  00185 Roma - Italy }\\
\emph{E-mail address:} \texttt{gambioli@mat.uniroma1.it}

\end{document}